\documentclass[11pt]{article}
\usepackage{amsthm}  
\usepackage{amsmath, amssymb, amsfonts}
\usepackage{graphicx}

\newtheorem{theorem}{Theorem}[section]

\newtheorem{definition}[theorem]{Definition}
\newtheorem{proposition}[theorem]{Proposition}
\newtheorem{corollary}[theorem]{Corollary}

\newtheorem{example}[theorem]{Example}

\renewcommand{\c}{{\boldsymbol{c}}}
\newcommand{\x}{{\boldsymbol{x}}}
\newcommand{\av}{{\boldsymbol{a}}}
\newcommand{\bv}{{\boldsymbol{b}}}

\newcommand{\C}{\mathcal{C}}
\newcommand{\F}{\mathcal{F}}

\renewcommand{\L}{\mathcal{L}}
\renewcommand{\P}{\mathcal{P}}
\newcommand{\R}{\mathbb{R}}

\newcommand{\N}{\mathcal{N}}

\newcommand{\ex}{\textrm{ex}}

\newcommand{\cl}{\tau}

\newcommand{\map}{\phi}

\newcommand{\comment}[1]{}

\begin{document}
\title{Pruning Processes and a New Characterization of \\Convex Geometries}
\author{Federico Ardila\thanks{Dept. of Mathematics, San
Francisco State University, San Francisco, CA, USA. 
({\tt federico@math.sfsu.edu}).}
\and
Elitza Maneva\thanks{IBM Almaden Research Center, 
San Jos\'e, CA, USA. ({\tt enmaneva@us.ibm.com}).}}
\date{}
\maketitle

\begin{abstract}
We provide a new characterization of convex geometries via a
multivariate version of an identity that was originally proved 
by Maneva, Mossel and Wainwright
for certain combinatorial objects arising
in the context of the $k$-SAT problem. We thus highlight the
connection between various characterizations of convex geometries and
a family of removal processes studied in the literature on random
structures.
\end{abstract}

\section{Introduction}

This article studies a general class of procedures in which the elements
of a set are removed one at a time according to a given rule. We
refer to such a procedure as a \emph{removal process}.  If every
element which is removable at some stage of the process remains
removable at any later stage, we call this a
\emph{pruning process}.
The subsets that one can reach through a pruning process have the
elegant combinatorial structure of a \emph{convex geometry}.  Our
first goal is to highlight the role of convex geometries in the
literature on random structures, where many pruning processes have
been studied without exploiting their connection to these objects.
Our second contribution is a proof that a generalization of a
polynomial identity, first obtained for a specific removal process in
\cite{MMW07}, provides a new characterization of pruning processes and
of convex geometries.  To prove this result we also show how a convex
geometry is equivalent to a particular kind of interval partition of
the Boolean lattice.

Two equivalent families of combinatorial objects, known as \emph{convex
geometries} and \emph{antimatroids}, were defined in the 1980s 
\cite{Edelman, KL84}.
The fact that these objects can be characterized via 
pruning processes has been known since then.  Some examples of
pruning processes considered at that time are the removal of vertices
of the convex hull of a set of points in $\R^n$, the removal of the
leaves of a tree, and the removal of minimal elements of a poset.
More recently various pruning processes have been studied in the
literature on random structures, and referred to also as peeling,
stripping, whitening, coarsening, identifying, etc.  A
typical example is the removal of vertices of degree less than $k$ in
the process of finding the $k$-core of a random (hyper)graph.

In \cite{MMW07}, a surprising identity was proved to hold for a
particular removal process which arises in the context of the $k$-SAT
problem. In this paper, we answer the question posed by Mossel
\cite{Mossel?} of characterizing the combinatorial structures that
satisfy (the multivariate version of) that identity: they are
precisely the convex geometries or equivalently the pruning processes. 
That is the content of our main result, Theorem
\ref{thm:partition} and Corollary \ref{cor:process}. It says that
any pruning process has the following two properties, and that in
fact either of these two properties characterizes pruning processes
among removal processes. 

\begin{itemize}
\item Suppose there is a subset $S$ of elements that we do not wish to remove.
Then there is a unique minimal set $\tau(S)$ achievable 
by the pruning process which contains $S$. 
\item Suppose to each element $e$ corresponds a weight $p_e$. For 
every set $S$ reachable by the pruning process, whose set of removable
elements is $U\subseteq S$, define the weight of $S$   
to be $\prod_{e\notin S} p_e \prod_{e \in U} (1-p_e)$. 
Then the sum of the weights of all reachable sets is 1.
\end{itemize}

Equivalently, after appropriate rewording, either of these properties
characterize convex geometries among set systems.



\medskip
\noindent {\bf Outline.} 
In Section \ref{sec:defn} we define precisely the equivalent concepts
of convex geometry, antimatroid, and pruning process.  We describe a
few different ways of looking at these objects: as set systems, as
languages generated by a set of circuits or pruning rules, and as
lattices. These different representations bring forward the
differences between various pruning processes that have been
studied. Section \ref{sec:identity} is devoted to our new
characterization of convex geometries.  Finally, in Section
\ref{sec:sat} we apply our results to the $k$-SAT problem and the
distribution considered in \cite{MMW07}, thereby generalizing Theorem
6 of \cite{MMW07}.

\medskip
\noindent {\bf Related work.} 
Antimatroids and convex geometries were first identified in the
context of lattice theory by Dilworth \cite{Dil40}. Since then they
have appeared in a variety of combinatorial situations. Two
particularly important treatments are Edelman and Jamison's convexity
approach \cite{Edelman}, and Korte and Lov\'asz's greedoid approach
\cite{KL84}.
Two good introductions to the subject are \cite{Bjorner} and \cite{Korte}.

Two notable examples of pruning processes on random structures
appear in the analysis of identifiable vertices in random hypergraphs
\cite{DN05}, and of $k$-cores in random hypergraphs
\cite{Mol04,PSW96}.  Additionally, in the analysis of satisfiability
problems such processes have appeared repeatedly --- most notably for
the pure-literal rule algorithm for $k$-SAT
\cite{BFU93,Mit97,Mol04,RPF99}, and the study of clustering of
solutions for XOR-SAT \cite{MRZ03} and $k$-SAT \cite{AR06,
MS07}. 

Pruning processes appear also in practical applications; for example in error correcting codes such as LDPC codes \cite{Gallager, LMSSS97} and LT codes \cite{Luby02} over the erasure 
channel. A unified analysis of the pruning processes in error-correcting codes
and the pure-literal rule is provided in \cite{LMS98}. 

Our work and in particular the implication $1 \rightarrow 3$ of our
main Theorem \ref {thm:partition} is related to previous work
of Aivaliotis, Gordon, and Graveman \cite{AGG01} and Gordon
\cite{Gor04}. For more information on this connection, see Section
\ref{sec:prob}.

\medskip
\noindent {\bf A word on terminology.} 
The objects of this paper have been studied under several different names. In particular, other authors have referred to pruning processes as \emph{shelling processes}. We prefer to avoid this name, which may lead to confusion with the other, better established notion of shelling in combinatorics. The term ``pruning" is more accurate, since a pruning process is equivalent to the successive removal of outermost elements of a convex geometry. A good example to keep in mind throughout the paper is the process of pruning of a tree by successively removing its leaves.

\section{Convex geometries and antimatroids}
\label{sec:defn}

Convex geometries and antimatroids are equivalent families of
combinatorial objects. Convex geometries provide a combinatorial
abstraction of the notion of convexity. Antimatroids describe pruning
processes, where we remove elements from (or add elements to) a set
one at a time, and once an element becomes available for removal, it
remains available until it is removed. There are many equivalent
definitions of these objects and a vast underlying theory
\cite{Bjorner, Korte}. We now present four points of view which we
will use.

\subsection{Convex sets and closure operations}

A \emph{convex geometry} is a pair $(E, \N)$ where $E$ is a set and
$\N \subseteq 2^E$ is a collection of subsets of $E$ satisfying:

\noindent{\bf (N1)} 
$E \in \N$.

\noindent{\bf (N2)} If $A,B \in \N$ then $A \cap B \in \N$.

\noindent{\bf (N3)} For every $A \in \N$ with $A \neq E$ there is an 
$x \notin A$ such that $A \cup x \in \N$.

The sets in $\N$ are called \emph{closed} or \emph{convex}.  
It is sometimes convenient to think of $\N$ as a poset ordered
by containment; this is a lattice.
We can then think
of (N3) as a property of \emph{accessibility from the top}: every
closed set can be obtained from $E$ by removing one element at a time
in such a way that every intermediate set in the process is also
closed. 

The \emph{closure} of $A\subseteq E$ is defined to be 
\[
\cl(A) = \bigcap_{\{C \in \N\, : \, C \supseteq A\}} C,
\]
which is the minimum closed set containing $A$.
It is easy to see that $\cl$ is, in fact, a \emph{closure operator};
that is, for all $A$ we have $A \subseteq \cl(A)$ and $\cl(\cl(A)) =
\cl(A)$, and for all $A \subseteq B$ we have $\cl(A) \subseteq
\cl(B)$. Also, a set $A$ is closed if and only if $\cl(A)=A$.

\begin{figure}[h]
\begin{center}
\begin{tabular}{ccc}
\includegraphics[width=1.5in]{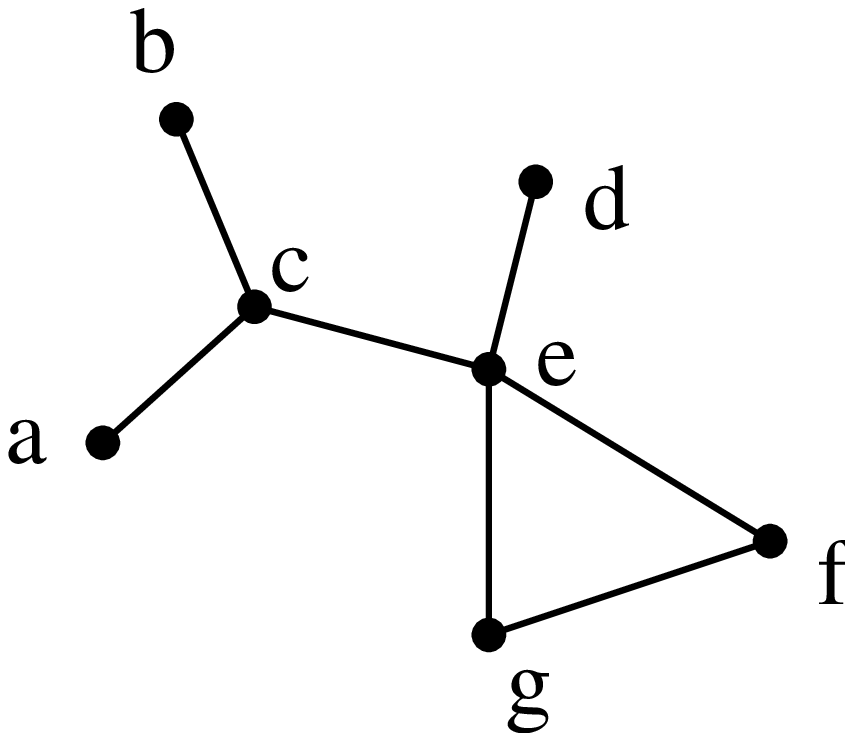} 
&~~~~~~~~~~~~~~~~~~& 
\includegraphics[width=2in]{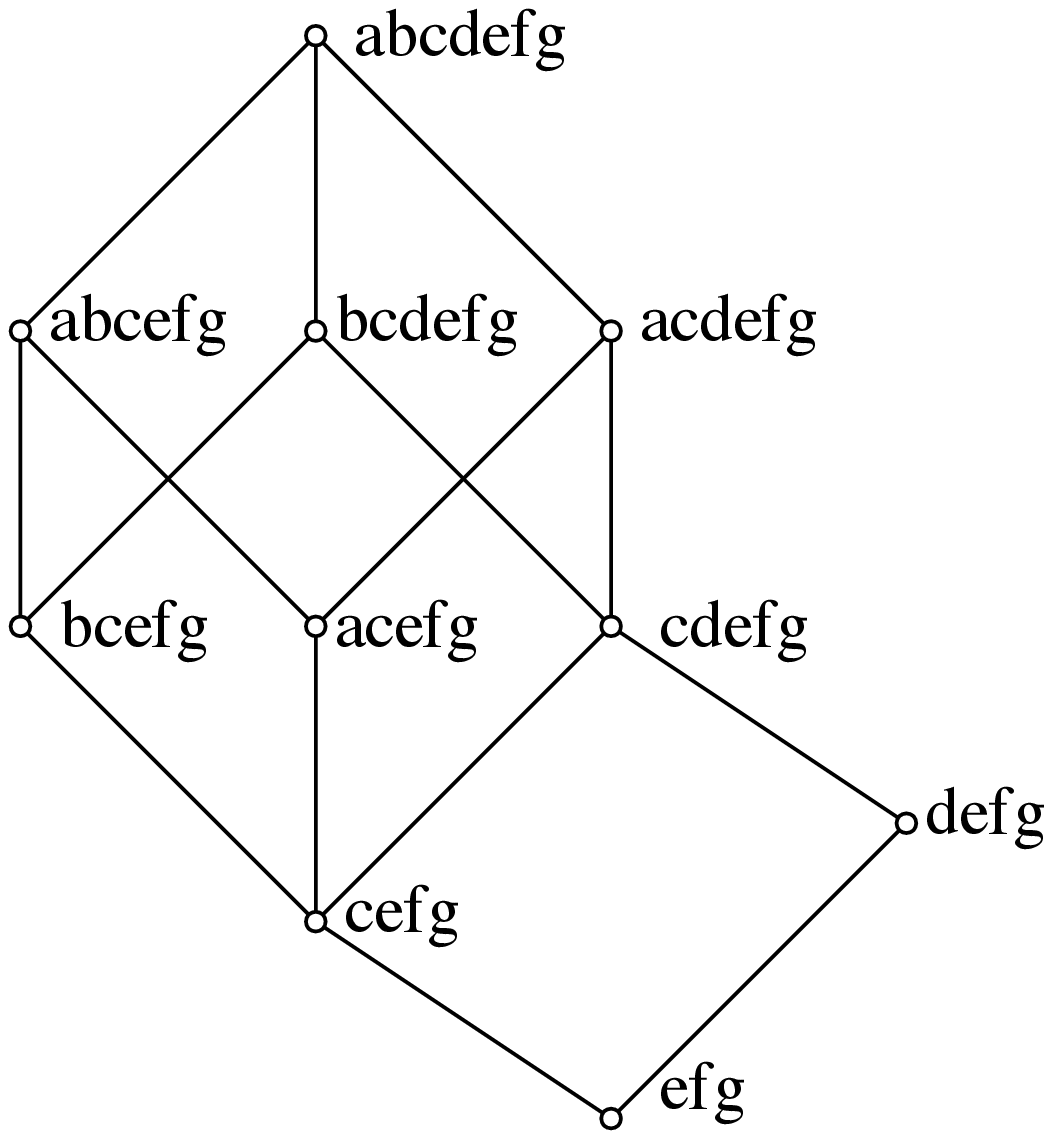} \\
\end{tabular}
\caption{\label{fig:example}
The set of configurations reachable by the process that 
successively removes leaves from a graph. The bottom configuration is the 2-core of the graph. }
\end{center}
\end{figure}

\begin{example}\label{example}
For a given graph, consider the subgraphs that can be 
obtained by successively removing leaves (nodes of degree 1). The vertex sets of these subgraphs are
the closed sets of a convex geometry. The minimal closed set is the 2-core of the graph. Figure \ref{fig:example} shows a specific graph and its lattice of closed sets; for example, the closure of the set $\{a,f,g\}$ is the set $\{a,c,e,f,g\}$. 
\end{example}

\comment{
\begin{example}\label{ex:conv}
Let $E$ be a set of points in $\R^n$. Say a subset $A \subseteq E$ is
\emph{closed} if it contains all the points $p \in E$ which are in the
convex hull of $A$. Then $(E, \N)$ is a convex geometry. The closure
of an arbitrary subset $A \subseteq E$ is the set of points of $E$
contained in the convex hull of $A$.
\end{example}

For the first configuration of Figure \ref{fig:config}, the only
non-closed sets are $\{c,d\}$ and $\{a,c,d\}$, whose closures are
$\{b,c,d\}$ and $\{a,b,c,d\}$, respectively.  For the second
configuration, the non-closed sets are $\{1,2\}, \{1,4\},
\{3,4\},\{1,2,4\},$ and $ \{1,3,4\}$.

\begin{figure}[h]
\begin{center}
\includegraphics[width=3.2in]{configs.eps} 
\caption{\label{fig:configs}
Two point configurations.}
\end{center}
\end{figure}
}

An \emph{extreme point} of a set $A$ is an element $a \in A$ which is
not in the closure of $A-a$. The set $\ex(A)$ of extreme points of $A$
is the unique minimal set whose closure is $A$.

\subsection{Antimatroids and pruning processes}\label{sec:antimatroids}

Next we define antimatroids, which are equivalent to convex geometries.
Let $E$ be a set, whose elements we regard as letters. A word over the
alphabet $E$ is called \emph{simple} if it contains no repeated
letters; let $E_s^*$ be the set of simple words over $E$.

An \emph{antimatroid} is a pair $(E, \L)$ where $E$ is a set and $\L
\subseteq E_s^*$ is a set of simple words satisfying:

\noindent{\bf (L1)}
If $\alpha\beta \in \L$ then $\alpha \in \L$; that is, any beginning
section of a word of $\L$ is in $\L$.

\noindent{\bf (L2)} If $\alpha, \beta \in \L$ and $|\alpha| > |\beta|$, then $\alpha$ contains a letter $x$ such that $\beta x \in \L$.

\noindent{\bf (L3)} If $\alpha, \beta \in \L$ and $x \in E$ are such that $\alpha x, \alpha\beta \in \L$ and $x \notin \beta$, then  $\alpha \beta x \in \L$.

Axiom (L1) says that $\L$ is left hereditary, (L2) is an exchange
axiom, and (L3) states that, as we build up a word of $\L$ from left
to right, any letter which can be added to the word at a certain
stage can still be added at any later stage.

The supports of the words in $\L$ are called the \emph{feasible}
subsets of $E$. The feasible subsets determine $\L$: a word is 
in $\L$ if and only if every initial segment of it is a feasible set.
The following theorem provides a
one-to-one correspondence between antimatroids and convex geometries.

\begin{theorem} \cite[Theorem III.1.3]{Korte}
Let $E$ be a finite set and $\F$ be a collection of subsets of
$E$. Then $\F$ is the collection of feasible sets of an antimatroid if
and only if $\F^* = \{E-F : F \in \F\}$ is the collection of closed
sets of a convex geometry.
\end{theorem}

By the above theorem, the feasible sets of the antimatroid
corresponding to our example of a convex geometry can be read on the
descending paths from the top of the lattice of closed sets. A letter
corresponds to every edge - this is the element that is removed. Any
set of removed elements is the complement of a convex set, thus it is
a feasible set. We will see below that the words of the antimatroid
can also be read on the descending paths.

\medskip

An alternative characterization of antimatroids starts by defining a set $H(x)
\subseteq{2^{E-x}}$ for every $x\in E$, which is a collection of 
\emph{alternative precedences}
for $x$; each precedence is a set not containing $x$. 
Let $\L$ be the
set of words on the alphabet $E$ such that $x$ can only appear after
at least one of its precedences has appeared:
\[
\L = \{x_1 \ldots x_k : \textrm{ for all } i \textrm{ there is a set } A \in H(x_i) \textrm{ with }A \subseteq \{x_1, \ldots, x_{i-1}\}.\}
\]

In the example of Figure \ref{fig:example} the alternative precedences
for $c$ are $\{a,b\}$, $\{a,e\}$, and $\{b,e\}$.  In general, for the
process of removing leaves to obtain the $2$-core of a graph, a vertex
of degree $d$ becomes removable as soon as $d-1$ of its neighbors are
removed. Thus its precedences are all the $(d-1)$-subsets of its set
of $d$ neighbors.

%






%



\medskip

A \emph{removal process} is a procedure in which the elements
of a set are removed one at a time according to a given rule. If every
element which is removable at some stage of the process remains
removable at any later stage, we call this a \emph{pruning process}.

We are now in a position to explain the correspondence between
pruning processes and antimatroids. Given an antimatroid $\L$ on
ground set $E$, we can consider each word $w$ of $\L$ as a removal
sequence, which instructs us to remove the elements of $w$ from left
to right. These removal sequences describe a pruning process: if an
element $x$ is removable at a certain stage described by word $w \in
\L$ (that is, if $wx \in \L$), then one of the alternative precedences
of $x$ appears in $w$. At any later stage the removal sequence $w'$
will contain $w$ as a prefix, and therefore will contain that
alternative precedence for $x$ as well.

Conversely, suppose we are given a pruning process on a set $E$. For
each element $x$ let the alternative precedences of $x$ be the subsets
$A \subseteq E$ such that $x$ is removable in $E-A$. Clearly the
antimatroid determined by these sets of alternative precedences
consists of the removal sequences in our pruning process.

Example \ref{example} gives rise to an antimatroid whose words are the
\emph{pruning sequences} of leaves that one can successively remove
from the graph. These correspond to the descending paths from the top
of the lattice of closed sets; the word indicates the elements that
are being removed as we walk down. For example the leftmost path in
the lattice of Figure \ref{fig:example} gives the word $dabc$, which
corresponds to a valid order of successively removing leaves from the
graph.

\subsection{Circuits and Paths}

The circuit description of a convex geometry is a good
way to reveal the pruning process which generates it.
A \emph{rooted set} is a set with a designated element called the root.
Any collection of rooted sets
gives rise to a convex geometry as follows. Suppose $\C$ is a family of rooted
subsets of $E$; label them $(A_i,a_i)$ where $a_i \in A_i$ and $A_i
\subseteq E$. Call a subset $S$ of $E$ \emph{full} if $A_i-a_i \subseteq S$
implies $a_i \in S$ for each $i$; that is, if the root of a set is
never the only element of the rooted set missing from $S$. Say a full
set $S$ is \emph{accessible from $E$} if there exists a sequence of
full subsets $E=S_0 \supset S_1 \supset \cdots \supset S_k=S$ with
$|S_i - S_{i+1}| = 1$ for each $i$. The following is, in a different
language, Lemma 3.7 of \cite{Korte}.

\begin{proposition} \cite{Korte} \label{prop:circuits}
Let $\C$ be a collection of rooted sets in $E$. The collection
$\N(\C)$ of full subsets of $E$ which are accessible from $E$ is the
collection of closed sets of a convex geometry on $E$.
\end{proposition}

The set of rooted sets can be interpreted as pruning rules. Let for
every $e\in E$, $\C_e=\{C\backslash\{e\}: (C, e) \in \C \})$. Then the
corresponding pruning process is the one in which an element $e$ is
removable if and only if at least one element has been removed from
each set in $\C_e$.

Conversely, from a convex geometry it is possible to recover a set of
rooted sets that generates it.
A \emph{free set} is one of the form $\ex(A)$.
A \emph{circuit} is a minimal set which is not free, and one can check
that each circuit $C$ has a unique element $a$ which is in the closure
of the remaining ones. This element is called the \emph{root} of the
circuit, and $(C,a)$ is called a \emph{rooted circuit}.

\comment{
In the first configuration of Figure \ref{fig:config}, the only sets
which are not free are $\{b,c,d\}$ and $\{a,b,c,d\}$, so the only
rooted circuit is $(\{b,c,d\},b)$. In the second configuration, the
rooted circuits are $(\{1,2,3\},3),(\{1,2,4\},2),(\{1,3,4\},3),$ and
$(\{2,3,4\},2)$.
}

\comment{
Note that in general, the collection of circuits of the convex
geometry $(E,\N(\C))$ is \textbf{not} equal to the collection $\C$ of
rooted sets that generated the convex geometry. However, these
\textbf{are} equal if we start with the collection $\C$ of circuits of
a convex geometry $(E, \N)$; it is not hard to see that in this case
we obtain $\N(\C)=\N$.
It is worth noting that for every convex geometry there exists a
minimal collection of rooted sets which generates it; these rooted 
sets are called its \emph{critical circuits}.
}

\comment{
There is a dual yet equivalent way of defining a convex geometry
the collection of full sets is the collection of closed sets of a
convex geometry, and every convex geometry arises in this way. 
}

Even though we will not need this fact, let us point out that the
collection of full sets of a family of rooted sets has a nice structure.

\begin{proposition}
\label{prop:fullsets} 
A collection $\F \subseteq 2^E$ is the collection of full sets
determined by a family of rooted sets if and only if it contains $E$
and is closed under intersection.
\end{proposition}
\begin{proof}
First, it is easy to see that the collection of full sets determined by a
family of rooted sets $\C$ contains $E$ and is closed under intersection.

Next, suppose $\F$ is a collection of subsets of $E$ that contains $E$
and is closed under intersection. We start with $\C$ being the
complete set of rooted sets on $E$. For every $F\in \F$ remove from
$\C$ all rooted sets $(A\cup a, a)$, where $A\subseteq F$ and $a\notin
F$.  We claim that $\F$ is the collection of full sets of $\C$.

It is immediate by the construction that every $F \in \F$ is a full set. 
It remains to show that there are no other full sets.
Suppose $D\subseteq E$ is a full set for $\C$. That means that
all rooted sets $(A\cup a, a)$ with $A\subseteq D$ and $a\notin D$ have
been removed. In particular $(D \cup a, a)$ has been removed for all $a
\notin D$. This implies that for every $a\notin D$ there exists $F_a\in
\F$ such that $D\subseteq F_a$ and $a\notin F_a$. Since $D = \cap F_a$ and
$\F$ is closed under intersection, $D$ is in $\F$.
\end{proof}

So far in this paper, rooted sets have played the role of 
circuits in a convex geometry. 
It is worth pointing out, however, that one can consider rooted sets as
\emph{paths} which generate a convex geometry in a different way, as follows.
Suppose $\P$ is a family of rooted subsets
of $E$ which we now call \emph{paths}; label them $(P_i,p_i)$ where $p_i \in P_i$ and $P_i
\subseteq E$. Let a subset $S$ of $E$ be \emph{path-full} if for every
$e \in S$ there exists $(P, e) \in \P$ such that $P\subseteq S$. Let it be \emph{path-closed}
if it is path-full, and accessible from $E$ by a sequence $E=S_0 \supset S_1 \supset \cdots \supset S_k=S$ of path-full subsets with $|S_i-S_{i+1}|=1$ for each $i$. Then the path-closed sets are the closed sets of a convex geometry, and every convex geometry arises in this way from a set of paths.

 In this context, the rooted sets can again be
interpreted as pruning rules. For every $e\in E$ let
$\P_e=\{P\backslash\{e\}: (P, e) \in \P \}$. Then the corresponding
pruning process is the one in which an element $e$ is removable if
and only if there is a set $P \in
\P_e$ such that every element of $P$ has been removed.

\comment{
Conversely, from a convex geometry it is possible to recover a set of
rooted sets that generates it in this way. An end point of a feasible
set $A$ of an antimatroid is an element $e\in P\backslash A$ such that
$A\cup \{e\}$ is also a feasible set. A \emph{path}, as defined in
\cite{Korte}, is a feasible set that has a single end point, or
equivalently it is a minimal feasible set containing $e$. A path with
an endpoint $e$ is considered as a rooted set $(P, e)$.
}

The pruning processes in the literature on random structures are
generated by rules for removing elements which can usually be
represented in a natural way through circuits or paths. While both
points of view are equivalent, sometimes one is more natural than
the other. For example, in finding the $k$-core of a graph \cite{PSW96}, 
a vertex becomes removable when at least one element has been 
removed from every $k$-subset of its neighbors (circuit rule). On the
other hand, in the case of identifiable vertices in hypergraphs
\cite{DN05} a vertex is removable if, in at least one of the
hyperedges in which it appears, every other vertex has been removed 
(path rule).

In the rest of the paper, all the rooted sets that appear play the role of circuits. 

\subsection{Lattices}

Finally, we can also think of convex geometries as meet-distributive
lattices. 
A lattice $L$ is \emph{meet-distributive} if for any element $x \neq
\widehat{0}$ the interval $[m(x),x]$ is a Boolean lattice, where $m(x)$
is the meet of the elements covered by $x$. Meet-distributive
lattices are precisely the posets of closed sets of convex geometries.
\cite[Prop. 8.7.5]{Bjorner} 

One might wonder whether something more specific can be said about the
convex geometries that arise from a set of circuits each of size at
most $k$. This question will be particularly natural in Section \ref{sec:sat},
where convex geometries are applied to the $k$-SAT problem for a fixed
value of $k$.

For $k=2$ the situation is very nice. Recall that a lattice $L$ is \emph{distributive} if 
$x \wedge (y \vee z) = (x \wedge y) \vee (x \wedge z)$ for any $x,y,z$ in
$L$.

\begin{proposition} \cite[Cor. 3.10]{Korte}
\label{prop:distributivity}
Let $\C$ be a set of rooted sets of size 2. The convex geometry $\N(\C)$
generated by these sets is a distributive
lattice. Conversely, every distributive lattice arises in this way.
\end{proposition}

For higher values of $k$, if we start with a collection $\C$ of rooted circuits of size $k$, the resulting
convex geometry $\N$ generally has additional rooted circuits of different
sizes. For the case of $k=2$, all circuits of the generated convex geometry have
size 2. However, for example, the convex geometry defined by the
rooted $3$-sets $(abc,b)$ and $(bde,d)$ also has $(acde,d)$ as a
circuit. 

If we have a bound on the size of the circuits of a convex geometry, we can make the following statement. A lattice $L$ is \emph{$k$-distributive} if 
$x \wedge (y_0 \vee \cdots \vee y_k) = (x \wedge y_0) \vee \cdots \vee (x \wedge y_k)$ for any $x,y_0, \ldots, y_k$ in $L$.

\begin{proposition} \cite[Cor. 4.3.]{Lib95}
\label{prop:kdistributivity}Let $k\geq 3$ be an integer.
If all circuits of a convex geometry have size at most $k$, then its poset of closed sets is a $(k-1)$-distributive lattice. Not every $(k-1)$-distributive lattice arises in this way.
\end{proposition}

\comment{
A generalization of one direction of Proposition
\ref{prop:distributivity} is proved in \cite{Lib95}. 
A lattice $L$ is \emph{$k$-distributive} if for any
$x, y_0, \dots, y_k \in L$, the following equation holds: $x \vee
\bigwedge_{i=0}^k y_i =
\bigwedge_{i=0}^k (x \wedge \bigvee_{j \neq i} y_j).$

\begin{proposition} \cite{Lib95}
\label{prop:n-distributivity}
Let $\C$ be a set of circuits of size at most $k$. The convex geometry
$\N(\C)$ is a $(k-1)$-distributive lattice. 
\end{proposition}
}

\section{Convex geometries as interval partitions of Boolean lattices}
\label{sec:identity}

In this section we describe our new characterization of convex geometries.
We show that convex geometries on a set $E$ are
characterized by the fact that they induce a certain partition of the
Boolean lattice $2^E$. This partition is encoded in a polynomial
identity which, as we will later see, generalizes Theorem
\ref{thm:polynomial} from \cite{MMW07}.

For any collection $\N \subseteq 2^E$ of subsets of $E$, and a set $A$
in $\N$, say that an element $a \in A$ is \emph{excludable from $A$}
if $A-a$ is in $\N$. Let $\ex(A)$ be the set of excludable elements of
$A$. When $\N$ is the collection of closed sets of a convex
geometry, $\ex(A)$ is the set of extreme points of $A$.

\begin{theorem} \label{thm:partition}
Let $\N \subseteq 2^E$ be a collection of subsets of a non-empty
set $E$. The following statements are equivalent:
\begin{enumerate}
\item $\N$ is the collection of closed sets of a convex geometry.
\item As $A$ ranges over $\N$ the intervals $[\ex(A), A]$ partition the
Boolean lattice $2^E$; that is, for every $D\subseteq E$ there
is a unique $A\in \N$ such that $\ex(A) \subseteq D \subseteq
A$. 
\item For any collection of $p_i$ and $q_i$ for $i\in E$ such that 
$p_i+q_i=1$ for all $i$, we have
\[
\sum_{A\in \N} ~~
\prod_{i\notin A} p_i \prod_{j \in \ex(A)} q_j=1.
\]
\end{enumerate} 
\end{theorem}

\begin{proof}
\noindent{\bf 1. implies 2.} 
Notice that if $A \in \N$ is such that $\ex(A) \subseteq D \subseteq
A$, then we have that  $A = \cl(\ex(A)) \subseteq \tau(D) \subseteq \cl(A) =
A$; so the only possible choice for $A$ is $A=\cl(D)$. It remains to
notice that, since $\ex(A)$ is the unique minimal set such that
$A=\cl(\ex(A))$, and $A = \cl(D)$, it follows that $D \supseteq
\ex(A)$, and therefore $D \in [\ex(A), A]$.

\noindent{\bf 2. implies 1.} We define the map $\map: 2^E \rightarrow \N$, as follows: for every
$D \subseteq E$, let $\map(D)$ be the unique element $A\in \N$ such
that $\ex(A) \subseteq D \subseteq A$. We need to show axioms
(N1)-(N3) of a convex geometry $(E,\N)$: $E$ is in $\N$, $\N$ is
closed under intersection, and every $A\in \N$ is accessible from $E$.

Axiom (N1) obviously holds, because $\map(E)=E$ is in
$\N$. To show (N3), we show that every set $A \in \N$ is accessible
from any superset $B \supseteq A$ that also belongs to $\N$.  It
suffices to prove that there exists an element of $B\backslash A$
that is excludable from $B$. If that were not the case, then
$\ex(B) \subseteq A$, and both of the intervals $[\ex(A),A]$ and $[\ex(B),B]$ 
would contain $A$, a contradiction.

Finally, we need to prove (N2), which states that $\N$ is closed under
intersection. First we prove the following statement:

\qquad \qquad \qquad \qquad
If $B \in \N$ and $A \subseteq B$ then $\map(A) \subseteq B$.

Suppose that we remove one element at a time from $B$ in any arbitrary
way, with the restriction that the intermediate sets in the process
must all be in $\N$ and contain $A$. We keep doing this until we
cannot continue anymore; suppose the set we obtain is $C$; by
construction, $B \supseteq C \supseteq A$. That means that every
element excludable from $C$ is in $A$, so $\ex(C) \subseteq A
\subseteq C$. Thus $C = \map(A)$ and we obtain the desired statement.

Now suppose $A_1$ and $A_2$ are in $\N$. From $A_1 \cap A_2 \subseteq
A_1$ we obtain that $\map(A_1 \cap A_2) \subseteq A_1$. Similarly
$\map(A_1 \cap A_2) \subseteq A_2$, so $\map(A_1 \cap A_2) \subseteq
A_1 \cap A_2$. But the reverse inclusion holds by definition, so we
must have equality. It follows that $A_1 \cap A_2$ is in $\N$.

\noindent{\bf 2. implies 3.} 
Observe that
$$ \sum_{D \subseteq E} ~
\prod_{i \not\in D} p_i \prod_{j \in D} q_j = 
\prod_{h \in E}(p_h + q_h)=1.$$
Therefore, since for every $D$ there is a unique $A$ such that $\ex(A)
\subseteq D \subseteq A$, it suffices to prove that for every $A \in
\N$: $$ \prod_{i\notin A} p_i \prod_{j \in \ex(A)} q_j=
\sum_{D \in [\ex(A), A]} ~
\prod_{i \not\in D} p_i \prod_{j \in D} q_j.  $$

This is easily seen to be true because:
\begin{eqnarray*}
\sum_{D \in [\ex(A),A]} ~
\prod_{i \not\in D} p_i \prod_{j \in D} q_j &=&
\prod_{i \not\in A} p_i \prod_{j \in \ex(A)} q_j 
\sum_{R \subseteq A\backslash \ex(A)} \left( 
\prod_{i \in A\backslash (\ex(A)\cup R)} p_i \prod_{j \in R} q_j \right) \\
&=& \prod_{i \not\in A} p_i \prod_{j \in \ex(A)} q_j
\prod_{h \in A\backslash\ex(A)}(p_h + q_h) \\
&=&
\prod_{i \not\in A} p_i \prod_{j \in \ex(A)} q_j.
\end{eqnarray*}

\noindent{\bf 3. implies 2.} Consider any set $D\subseteq E$,
and let $p_a=0$ if $a\in D$ and $p_a=1$ otherwise. The equality becomes:
\begin{eqnarray*}
1&=& \sum_{A\in \N}~~ \prod_{i \notin A} p_i \prod_{j \in \ex(A)} q_j
\\
&=&\sum_{A\in \N ~:~ \ex(A)\subseteq D \subseteq A} 1 
\end{eqnarray*}
Therefore there is exactly one set $A\in \N$ for which 
$\ex(A)\subseteq D \subseteq A$.
\end{proof}

The following corollary gives a
characterization of pruning processes among removal processes:


\begin{corollary}
\label{cor:process}
A removal process on a set $E$ is a pruning process if and only if,
for each subset $S$ of $E$, there is a unique minimal set $\cl(S)$
containing $S$ which is achievable by the removal process.
\end{corollary}
\begin{proof}
As outlined in Section \ref{sec:antimatroids}, a pruning process
gives rise to a convex geometry, and in that case $\cl(S)$ is just the
convex closure of $S$.  For the other direction, let $\N$ consist of
the sets achievable by the removal process; it suffices to show that
$\N$ is the collection of closed sets of a convex geometry. We will
show that property 2 of Theorem \ref{thm:partition} holds. For any set
$S\subseteq E$, it holds that $\ex(\cl(S)) \subseteq S \subseteq
\cl(S)$, because if there is an excludable element of $\cl(S)$ that is not in 
$S$, then $\cl(S)$ would not be the minimal set containing $S$.
Furthermore, any set $T\in \N$ for which $\ex(T)\subseteq
S\subseteq T$ is a minimal set containing $S$ because all of its excludable
elements are in $S$. Since there is a unique such set, it is $\cl(S)$.
\end{proof}

We conclude this section by offering a probabilistic interpretation of
property 3 of Theorem \ref{thm:partition}.

\subsection{A probabilistic interpretation}
\label{sec:prob}

Let $(E,\N)$ be a convex geometry, and fix $0 \leq p_e,q_e \leq 1$
with $p_e+q_e=1$ for each element $e$ of $E$.
Define a probability distribution $\pi_1$ on the subsets of $E$ by
independently deleting element $e$ with probability $p_e$ and keeping
it with probability $q_e$:
\[
\Pr_{\pi_1}(A) =  \prod_{i \notin A} p_i \prod_{j \in A} q_j, \qquad A
\subseteq E.
\]
Define a probability distribution $\pi_2$ on the convex sets of $E$ by:
\[
\Pr_{\pi_2}(A) =  \prod_{i \notin A} p_i \prod_{j \in \ex(A)} q_j,
\qquad A \in \N.
\]
The implication $1 \rightarrow 3$ of Theorem \ref{thm:partition} tells
us that $\pi_2$ is, indeed, a probability distribution. Furthermore,
in order to sample from $\pi_2$, it suffices to sample from $\pi_1$ and
compute the closure of the obtained set.

\begin{theorem}\label{th:prob}
Let $f(A)$ be a function defined on the subsets $A$ of $E$ which
depends only on $\cl(A)$.
Then the expected value of $f$ when we sample from the
distribution $\pi_1$ on all subsets of $E$, equals the expected value
of $f$ when we sample from the distribution $\pi_2$ on the convex
sets of $E$.
\end{theorem}

\begin{proof} Assuming that $f(D) = f(\cl(D))$, the identity
\begin{equation}\label{probs}
\sum_{D \subseteq E} ~~f(D)
\prod_{i\notin D} p_i \prod_{j \in D} q_j=
\sum_{A\in \N} ~~f(A)
\prod_{i\notin A} p_i \prod_{j \in \ex(A)} q_j.
\end{equation}
can be established in exactly the same way as implication $1
\rightarrow 3$ of Theorem \ref{thm:partition}.
\end{proof}

We note that Aivaliotis, Gordon, and Graveman \cite{AGG01} and Gordon
\cite{Gor04} studied the problem of choosing a random subset of a
convex geometry under the distribution $\pi_1$. They related the
expected rank of this random subset to the Tutte polynomial of the
antimatroid. In particular, they discovered (\ref{probs}) in a
special case which is no simpler than the general case.

The fact that $\pi_2$ is a probability distribution on $\N$ is not explicitly stated in \cite{AGG01} or \cite{Gor04}, and neither is the probabilistic interpretation of the right hand side of (\ref{probs}).
However, these two results follow very easily from that work.
Our theorem that the probabilistic property of Theorem \ref{th:prob}
(or the weaker condition 3 of Theorem \ref{thm:partition})
characterizes convex geometries is new.


\section{Convex geometries in the k-SAT problem}
\label{sec:sat}

Let $F$ be a Boolean formula such as
\[
F=(\bar{x}_1 \vee \bar{x}_2 \vee x_3) \wedge (x_2 \vee \bar{x}_3 \vee \bar{x}_4).
\]
We can assume that $F$ is written in \emph{conjunctive normal form} as
a conjunction of certain clauses $C$ in the variables $V$ and their
negations. The \emph{Boolean satisfiability problem (SAT)} is to
determine whether there is some assignment of TRUE $(1)$ and FALSE
$(0)$ to the variables which makes the entire formula true. The
\emph{$k$-SAT problem} is the same problem when restricted to
formulas with clauses of a fixed size $k$.  For $k=2$ there is a
polynomial time algorithm for deciding satisfiability, however for
$k\ge 3$ the problem is NP-complete.

In their analysis of the Survey Propagation algorithm
\cite{MPZ02,BMZ05} for 3-SAT, Maneva et al \cite{MMW07} discovered a
polynomial identity that holds for any SAT problem and any satisfying
assignment. To define this identity first we need to introduce the
concept of \emph{partial assignments}, where to each variable is
assigned one of the values $0,1,$ or $*$; the value $*$ indicates that
a variable is unassigned and free to take either value.  Say that a
partial assignment $\x$ is
\emph{invalid} for a clause $C$ if plugging $\x$ into $C$ gives either
$0 \vee 0 \vee \cdots \vee 0$ (which makes the clause invalid) or $0
\vee \cdots \vee 0 \vee * \vee 0 \vee \cdots \vee 0$ (where the $*$ is not free
to take either value). A partial assignment $\x$ is \emph{valid} for a
formula if it is valid for all its clauses.  For example, some valid
partial assignments for the formula $F=(\bar{x}_1 \vee \bar{x}_2 \vee
x_3) \wedge (x_2 \vee \bar{x}_3 \vee \bar{x}_4)$ are $(1,1,1,1),
(*,1,*,*)$ and $(1,*,*,1)$, and some invalid partial assignments are
$(1,1,0,*)$ and $(*,*,1,1)$.

\begin{definition} 
\label{def:poset}
Given a Boolean formula $F$, the poset $P(F)$ of valid partial
assignments is defined by decreeing that $\av$ covers $\bv$ if $\bv$ is
obtained from $\av$ by switching a $0$ or $1$ to a $*$.
\end{definition}

\begin{figure}[h]
\begin{center}
\includegraphics[width=4.4in]{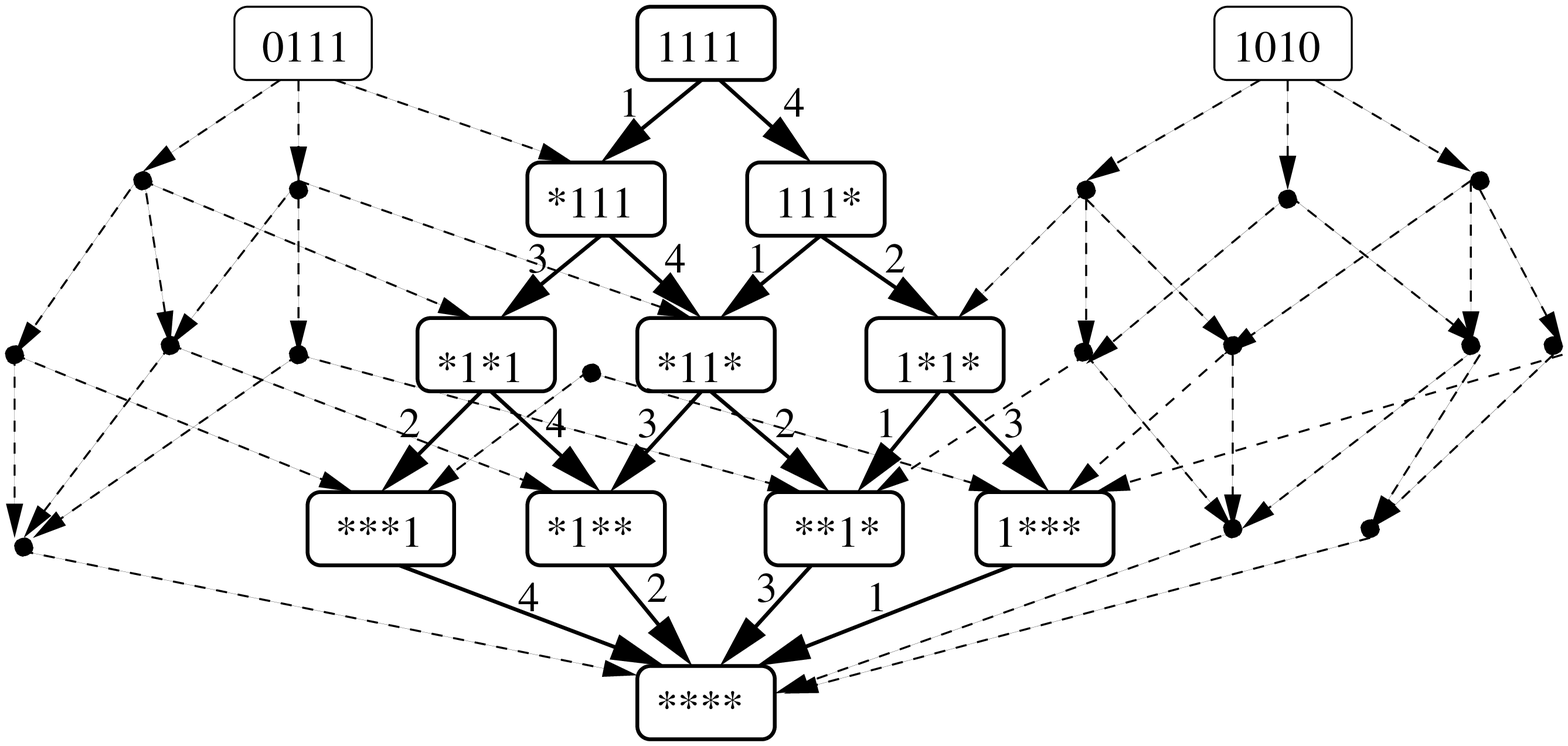} 
\caption{\label{fig:lattice}
Some of the valid partial assignments for the 
formula $(\bar{x}_1 \vee \bar{x}_2 \vee x_3) \wedge ( x_2 \vee
\bar{x}_3 \vee \bar{x}_4)$.  Highlighted are the assignments below the
satisfying assignment $(1,1,1,1)$. Edges are labeled with
the index of the variable whose value differs in the two adjacent
assignments.}
\end{center}
\end{figure}

Figure \ref{fig:lattice} shows part of the poset of valid partial
assignments for the formula $F$ above. Note that in this example, somewhat surprisingly, $(1,1,1,1)$ is not
greater than $(1,*,*,1)$ in $P(F)$, because to stay valid one must
switch $x_2$ and $x_3$ from $1$ to $*$ simultaneously: $(1,*,1,1)$ and
$(1,1,*,1)$ are invalid.

Definition \ref{def:poset} suggests that, given a valid partial
assignment $\av$, we call a coordinate $i$ either:

\noindent (a) a \emph{star} $*$,

\noindent (b) \emph{unconstrained} if $a_i\in\{0,1\}$ and setting 
$a_i=*$ keeps the assignment valid, or

\noindent (c) \emph{constrained} if $a_i\in\{0,1\}$ and setting $a_i=*$ 
gives an invalid assignment; that is, if $a_i$ is the only satisfying
variable in some clause of $F$. 

\noindent 
Let $S(\av)$, $U(\av)$, 
$C(\av)$, and $N(\av)$ be the sets of star, unconstrained, constrained, and numerical
variables of $\av$, respectively; so $V=S(\av)
\cup N(\av)$ and $N(\av)=U(\av) \cup C(\av)$.

Maneva et al \cite{MMW07} defined the \emph{weight} of a partial
assignment $\av$ to be
\[
W(\av) = p^{|S(\av)|}q^{|U(\av)|},
\]
where $p$ and $q$ are parameters in the interval $[0,1]$. They
considered the probability distribution which assigns to $\av$ a
probability proportional to $W(\av)$ for every valid partial
assignment $\av$.  The survey propagation algorithm was then proved to
be equivalent to applying the belief propagation marginalization
heuristic \cite{Pearl88} to this distribution with suitably chosen $p$
and $q$. This distribution has the following property, which should
not look surprising in view of Theorem \ref{thm:partition}:

\begin{theorem} \cite{MMW07} 
\label{thm:polynomial}
For any satisfying assignment $\av$ of a Boolean formula $F$ and
$p+q=1$,
\[
\sum_{\bv \leq \av} p^{|S(\bv)|} q^{|U(\bv)|} =  1, 
\]
summing over all valid partial assignments $\bv$ which are less than $\av$
in $P(F)$; that is, summing over the subposet $P(F)_{\leq \av}$.
\end{theorem}

Thus the probability distribution on partial assignments with $p>0$
may be regarded as a ``smoother'' version of the uniform distribution
over satisfying assignments, which corresponds to the case $p=0$: if we choose
a valid partial assignment $\bv$ at random, then the probability of
being under $\av$ is the same for any satisfying assignment
$\av$. Another consequence of the above theorem is that if the total
weight of all valid partial assignments is less than 1, then the formula has no
satisfying assignment.
In recent work of Sinclair and the second author \cite{MS07}, this fact
was used in conjunction with the first-moment
method to bound the probability of satisfiability of a random
SAT formula with clauses of sizes 2 and 3.

\medskip

Consider the following experiment:

1. in a valid assignment $\av$, change a random unconstrained variable to $*$, and 

2. repeat until there are no unconstrained variables.\\
This procedure has been referred to as ``peeling'', ``whitening'', 
``coarsening'' and
``pruning''. We now recognize it as a pruning process on the set
of variables which have numerical values in $\av$. At each stage,
we are allowed to remove an unconstrained variable; notice that if a 
variable becomes 
unconstrained, it remains unconstrained throughout this process.

This experiment is equivalent to taking a random path from $\x$ down the
partial order $P(F)$, by choosing at each step a random partial
assignment that is covered by the current one. For a fixed choice of
$\av$, any such path terminates at the same partial assignment, which
is known as a ``core''. (Note, however, that different $\av$ may lead
to different core assignments.)

Achlioptas and Ricci-Tersenghi \cite{AR06} examined the above pruning process and proved that, for $k \ge 9$ and a formula chosen from a particular distribution of interest, there is a high probability that the  process will terminate before removing all variables.
This is not known to hold for $k=3$.

With the above description of the removal process, our next result follows easily.

\comment{
Let $F$ be a SAT formula with variables $V$, and let $\av$ be a valid
(possibly partial) assignment for $F$. For each valid partial
assignment $\bv$ which is less than $\av$ in the poset $P(F)$, let the
sets of star, numerical, unconstrained, and constrained variables of
$\bv$ be $S(\bv), N(\bv), U(\bv), C(\bv)$, respectively; so $V=S(\bv)
\cup N(\bv)$ and $N(\bv)=U(\bv) \cup C(\bv)$.
}

\begin{theorem} \label{thm:sat-antimatroid}
Let $F$ be a SAT formula with variables $V$, and let $\av$ be a
valid (possibly partial) assignment for $F$. Let
\[
\N = \{N(\bv) \, : \, \bv \textrm{ is a valid partial assignment such that } \bv\leq \av\}.
\]
Then $(N(\av),\N)$ is a convex geometry.
%
Conversely, every convex geometry arises in this way from a valid
assignment for a SAT formula.
\end{theorem}

\begin{proof}
We show that this statement is equivalent to Proposition \ref{prop:circuits}.
Consider the clauses of $F$ with a unique satisfying variable in $\av$,
which give $0 \vee \cdots\vee 0 \vee 1 \vee 0
\vee \cdots \vee 0$ when we plug $\av$ into them. If $C$ is the set of variables in such a clause
(which must be a subset of $N(\av)$) and $v$ is the unique satisfying
variable, form a rooted set $(C,v)$. Then $(N(\av),\N)$ is clearly the
convex geometry generated by these rooted sets. Conversely, given a
convex geometry, one can encode its rooted sets into the clauses of a
SAT formula with a valid assignment.
\end{proof}

The convex geometry corresponding to assignment $(1,1,1,1)$ in Figure
\ref{fig:lattice} is the collection of sets of assigned variables in assignments
lying below $(1,1,1,1)$: 
$$\{
\{1,2,3,4\},  
\{2,3,4\},  
\{1,2,3\},  
\{2,4\},  \{2,3\},  \{1,3\},  
\{4\}, \{2\}, \{3\}, \{1\}, \emptyset\}$$
and the words of the corresponding antimatroid can be read out by going down the
directed edges; the feasible sets are:
$$\{ \emptyset, \{1\}, \{4\},   \{1,3\}, \{1,4\}, \{2,4\}, \{1,2,3\}, \{1,3,4\}, \{1,2,4\}, \{2,3,4\},  \{1,2,3,4\}\}$$

In the particular case of 2-SAT, the convex geometry is very special. By
Proposition \ref{prop:distributivity}, the poset $P(F)_{\leq a}$ of
valid partial assignments is a distributive lattice.

Notice that a SAT formula $F$ generally has several different
valid partial assignments, and each assignment $\av$ gives rise to
a convex geometry $G(F,\av)$.  These different convex geometries fit
together nicely, as seen in Figure \ref{fig:lattice}.  If $G(F,\av)$
and $G(F,\bv)$ have a non-empty intersection, then their intersection
is the convex geometry $G(F,\c)$ for the unique element $\c$ with
maximal $N(\c)$ for which $c_i=*$ if $a_i \neq b_i$, and $c_i=a_i=b_i$
otherwise. 

The machinery that we have built up now provides a more
illustrative multivariate version of Maneva, Mossel, and Wainwright's
Theorem \ref{thm:polynomial} on the probability distribution determined
by a SAT problem $F$ and a valid assignment $\av$.  
%
More importantly, in view of Theorem \ref{thm:partition}, it
tells us that the identity of Theorem \ref{thm:polynomial} holds
precisely \textbf{because} a SAT problem gives rise to a convex
geometry. Therefore convex geometries are really the context in which
this identity should be understood.

\begin{theorem} \label{thm:multivariate}
For a valid partial assignment $\bv$ of a Boolean formula $F$ with variables $V$, let
$S(\bv), U(\bv)$ and $C(\bv)$ denote the sets of star, unconstrained, and
constrained variables of $\bv$ in $F$, respectively. Let $p_i$ and $q_i$ be 
such that $p_i+q_i=1$ for all $i \in V$.
Then, for any valid assignment $\av$ of a Boolean formula
$F$,
\[
\sum_{\bv \leq \av}  \,\,\prod_{i \in S(\bv)} p_i \prod_{j \in U(\bv)} q_j = 1
\]
summing over all valid assignments $\bv$ which are less than $\av$
in $P(F)$.
\end{theorem}

\begin{proof}
The result follows directly from Theorems \ref{thm:partition} and \ref{thm:sat-antimatroid}.

\comment{
Let $(N(\av), \N)$ be the convex geometry corresponding to $F$ by
Theorem \ref{thm:sat-antimatroid}. We can apply the second statement of
Theorem \ref{thm:partition} to the convex geometry, thereby
decomposing $2^{N(\av)}$ into intervals $[\ex(A),A]$ for $A$
closed. Each closed set is of the form $N(\bv)$ for a valid partial
assignment $\bv \leq \av$. Now, $N(\bv) - i$ is closed if and only if
switching the $i$th coordinate of $\bv$ from $b_i$ to $*$ results in a
valid partial assignment; this occurs if and only if $i$ is
unconstrained.  Therefore $\ex(N(\bv)) = U(\bv)$. Applying the third
statement of Theorem \ref{thm:partition}  gives 
\[
\sum_{\bv \leq \av} ~~ \prod_{i\in S(\bv)
} p_i 
\prod_{j \in \ex(N(\bv))} q_j=1,
\]
which is the desired result.
}
\end{proof}

\section{Acknowledgments}

We would like to thank Elchanan Mossel for posing the question of
characterizing the combinatorial objects satisfying Theorem \ref{thm:polynomial}, 
and Laci Lov\'asz and Martin Wainwright for helpful discussions. We would
also like to thank the referee for very useful suggestions for improving
the exposition.

\bibliographystyle{plain}
\bibliography{pruning}

\end{document}